\documentclass[12pt]{amsart}
\usepackage{amssymb}
\usepackage[all]{xy}

\newcommand{\issuenumber}{11}
\newcommand{\issuemonth}{December}
\newcommand{\issueyear}{2004}

\newtheorem{thm}{Theorem}[section]
\newtheorem{prob}[thm]{Problem}
\newtheorem{cor}[thm]{Corollary}
\newtheorem{lem}[thm]{Lemma}

\newtheorem{issue}{Issue}

\theoremstyle{definition}

\theoremstyle{remark}

\newcommand{\ed}{\end{thebibliography}\general\end{document}}

\renewcommand{\b}{\mathfrak{b}}

\newcommand{\p}{\mathfrak{p}}

\newcommand{\NON}{{\mathsf   {NON}}}
\newcommand{\COF}{{\mathsf   {COF}}}

\newcommand{\dnannouncement}[1]{[\S\ref{#1} below]}

\newcommand{\M}{\mathcal{M}}

\newcommand{\cov}{\mathsf{cov}}

\newcommand{\fo}{\mathfrak{od}}

\renewcommand{\b}{\mathfrak{b}}
\renewcommand{\t}{\mathfrak{t}}

\renewcommand{\split}{\mathsf{Split}}
\newcommand{\bq}{\begin{quote}}
\newcommand{\eq}{\end{quote}}

\newcommand{\B}{\mathcal{B}}
\newcommand{\BG}{\B_\Gamma}

\newcommand{\sone}{\mathsf{S}_1}    \newcommand{\sfin}{\mathsf{S}_{fin}}
\newcommand{\ufin}{\mathsf{U}_{fin}}

\newcommand{\nin}{\not\in}

\newcommand{\NN}{{{}^{\naturals}\naturals}}
\newcommand{\naturals}{{\mathbb N}}
\newcommand{\N}{\naturals}

\newcommand{\by}[2]{\par\hfill\emph{#1}, #2}
\newcommand{\nby}[1]{\par\hfill\emph{#1}}
\newcommand{\Tau}{\mathrm{T}}
\newcommand{\CE}{\textsc{CE}}

\newcommand{\be}{\begin{enumerate}}
\newcommand{\ee}{\end{enumerate}}
\newcommand{\bi}{\begin{itemize}}
\newcommand{\ei}{\end{itemize}}

\newcommand{\general}{\small\vfill\par\noindent\hrulefill\par
\noindent\textbf{Previous issues.} The first issues of this
bulletin, which contain general information (first issue), basic
definitions, research announcements, and open problems (all
issues) are available online, on \arx{math.GN/$x$}, where $x$ is
\texttt{0301011}, \texttt{0302062}, \texttt{0303057},
\texttt{0304087}, \texttt{0305367}, \texttt{0312140},
\texttt{0401155}, \texttt{0403369}, \texttt{0406411},
and \texttt{0409072},
respectively, for issues number $1$ to $10$.\\[0.1cm]
\textbf{Contributions.}
Please submit your contributions (announcements, discussions, and open problems)
by e-mailing us. It is preferred to write them
in \LaTeX{}.
The authors are urged to use as standard notation as possible, or otherwise give
the definitions or a reference to where the notation is explained.
Contributions to this bulletin would not require any transfer of copyright,
and material presented here can be published elsewhere.\\[0.1cm]
\textbf{Subscription.}
To receive this bulletin (free) to your
e-mailbox, e-mail us:\\
{tsaban@math.huji.ac.il}
}

\newcommand{\nArxPaper}[5]{\subsection{#2}{#4}\par\hfill{\arx{#1}}\par\hfill\emph{#3}}

\newcommand{\nAMSPaper}[5]{\subsection{#2}{#4}\par\hfill{\texttt{#1}}\par\hfill\emph{#3}}
\newcommand{\SPMBul}{\textbf{$\mathcal{SPM}$ Bulletin}}

\newcommand{\arx}[1]{\texttt{http://arxiv.org/abs/#1}}
\newcommand{\url}[1]{\bq\texttt{#1}\eq}
\newcommand{\online}[1]{The paper is available online at \url{#1}}

\newcommand{\probmonth}{\emph{Problem of the month}}

\title[$\mathcal{SPM}$ Bulletin \textbf{\issuenumber} (\issuemonth{} \issueyear)]{%
$\mathcal{SPM}$ Bulletin\\[0.5cm]
Issue number \issuenumber: \issuemonth{} \issueyear{} \CE{}}

\begin{document}
\maketitle

\tableofcontents

\section{Editor's note}

This issue contains, in addition to the usual contents, a special festive
announcement: A book.
This book by Banakh and Zdomsky is announced in \dnannouncement{SFbook},
and seems to be the first in a planned series by these authors.
We believe that the book will become a cornerstone in many future mathematical
investigations, in particular in the field of infinite-combinatorial topology.
The book's preliminary version is available online, as seen in the announcement,
and the readers of the \SPMBul{} are encouraged to take a look and make comments.

Zdomsky has also made two detailed contributions to this issue.
This is the ideal form of a contribution to the \SPMBul{}, and
we urge all contributors to consider this possibility from time to
time.

Finally, the BEST 2005 Conference is getting close -- have a look at the
corresponding announcement.

\medskip

Contributions to the next issue are, as always, welcome.

\medskip

\by{Boaz Tsaban}{tsaban@math.huji.ac.il}

\hfill \texttt{http://www.cs.biu.ac.il/\~{}tsaban}

\section{Research announcements}

\nAMSPaper{http://www.ams.org/journal-getitem?pii=S0002-9939-04-07744-5}
{On subclasses of weak Asplund spaces}
{Ondrej F.\ K.\ Kalenda and Kenneth Kunen}
{Assuming the consistency of the existence of a measurable cardinal,
it is consistent to have two Banach spaces, $X,Y$,
where $X$ is a weak Asplund space such that $X^{*}$ (in the weak*
topology) in not in Stegall's class, whereas $Y^{*}$
is in Stegall's class but is not weak* fragmentable.
}

\nArxPaper{math.LO/0409110}
{The number of translates of a closed nowhere dense set required to cover a Polish group}
{Arnold W.\ Miller and Juris Steprans}
{For a Polish group $G$ let $\cov_G$ be the minimal number of translates of a fixed
closed nowhere dense subset of $G$ required to cover $G$. For many locally compact
$G$ this cardinal is known to be consistently larger than $\cov(meager)$ which is
the smallest cardinality of a covering of the real line by meagre sets. It is
shown that for several non-locally compact groups $\cov_G=\cov(meager)$. For
example the equality holds for the group of permutations of the integers, the
additive group of a separable Banach space with an unconditional basis and the
group of homeomorphisms of various compact spaces.}

\nAMSPaper{http://www.ams.org/journal-getitem?pii=S0002-9939-04-07685-3}
{More on convexity numbers of closed sets in $\mathbb{R}^n$}
{Stefan Geschke}
{The \emph{convexity number} of a set $S\subseteq\mathbb R^n$ is the least size of
a family $\mathcal F$ of convex sets with $\bigcup\mathcal F=S$.
$S$ is \emph{countably convex} if its convexity number is countable.
Otherwise $S$ is \emph{uncountably convex}.

Uncountably convex closed sets in  $\mathbb R^n$ have been studied recently
by
Geschke, Kubi\'s, Kojman and Schipperus.
Their line of research is continued in the present article.
We show that for all $n\geq 2$,
it is consistent that there is an uncountably convex closed
set $S\subseteq\mathbb R^{n+1}$ whose convexity number is strictly
smaller than all convexity numbers of uncountably convex subsets of
$\mathbb R^n$.

Moreover, we construct a closed set $S\subseteq\mathbb R^3$ whose convexity
number is $2^{\aleph_0}$ and that has no uncountable $k$-clique for any $k>1$.
Here $C\subseteq S$ is a \emph{$k$-clique} if the convex hull of no
$k$-element subset of $C$ is included in $S$.
Our example shows that the main result of the above-named authors, a closed
set $S\subseteq\mathbb R^2$ either has a perfect $3$-clique
or the convexity number of $S$ is $<2^{\aleph_0}$ in some
forcing extension of the universe, cannot be extended to higher dimensions.
}

\label{SFbook}
\nAMSPaper{\noindent \scriptsize http://www.franko.lviv.ua/faculty/mechmat/Departments/Topology/booksite.html}
{A new book: Coherence of Semifilters}
{Taras Banakh and Lubomyr Zdomsky}
{The book is devoted to studying the (sub)coherence relation on semifilters,
that is families of infinite subsets of $\N$, closed under taking almost supersets.
On the family of ultrafilters the coherence relation was introduced in eighties by Andreas Blass,
who formulated his famous principle, the Near Coherence of Filters (NCF),
that found many non-trivial applications in various fields of mathematics.

In the book the (sub)coherence relation is treated with help of cardinal functions defined on
the lattice SF of semifilters.
Endowed with the Lawson topology the lattice SF becomes a supercompact topological space.
It can be interesting for topological algebraists because any reasonable binary operation on natural
numbers induces a right-topological operation on SF in the same way as it does on the Stone-Cech
compactification $\beta\N$.

A preliminary version of the book is available online:
}

\section{Characterization of topological spaces with (strictly) $o$-bounded
free topological group}

This announcement is devoted to the problem
 of characterization of Tychonov spaces $X$ such that the corresponding free
(abelian) topological group $F(X)$ ($A(X)$) is [strictly] $o$-bounded posed in \cite{HRT},
see \cite{HRT} or \cite{Tk98} for corresponding definitions.
All topological spaces are assumed to be Tychonov.

The  characterization we present here
 involves the concepts of  selection principles
and the Menger game on a uniform space.
 A uniform space $(X,\mathcal U)$ is said to to have the property
$\bigcup_{fin}(\mathcal O, \mathcal X)$, where $\mathcal X\in\{\mathcal O,\Omega,\Gamma\}$,
if for every sequence $(u_n)_{n\in\omega}$ of open uniform covers
of $X$ there exists a sequence $(v_n)_{n\in\omega}$ such that $\{\cup v_n:n\in\omega\}\in\mathcal X$
and each $v_n$ is a finite subfamily of $u_n$.
Next,  the Menger game
 on a uniform space $(X,\mathcal U)$ is obtained from the Menger game on the
 corresponding topological space (see \cite{Sch} for its definition) by
restriction of the choice of the first player to open uniform covers of $X$.
For a topological (uniform) space $X$ ($(X,\mathcal U)$) we shall simply
write $II\uparrow M(X)$ ($II\uparrow M(X,\mathcal U)$) in place
of ``the second player has a winning strategy in the
Menger game on $X$ ($(X,\mathcal U)$).
Let us note, that a topological group $G$ is (strictly) $o$-bounded,
if the uniform space $(G,\mathcal V^\ast)$ has the property $\bigcup_{fin}(\mathcal O,\mathcal O)$
($II\uparrow M(G,\mathcal V^\ast)$), where $\mathcal V^\ast$ is the right uniformity of
$G$.

Let $X$ be a Tychonov space. Recall from \cite{En} that the  uniformity
$\mathcal U$ on $X$ is called \emph{universal}, if it generates the topology
of $X$ and contains all uniformities on $X$ with this property.
Throughout this section the universal uniformity of a topological space
$X$ will be denoted by $\mathcal U(X)$.
We are in a position now to present the above mentioned characterizations.

\begin{thm} \label{main}
Let $X$ be a Tychonov space. Then the following conditions are equivalent:
\begin{itemize}
\item[$(1)$]  $F(X)$ is (strictly) $o$-bounded;
\item[$(2)$] $F(X)^n$ is (strictly) $o$-bounded for all $n\in\mathbb N$;
\item[$(3)$]  $A(X)$ is (strictly) $o$-bounded;
\item[$(4)$] $A(X)^n$ is (strictly) $o$-bounded for all $n\in\mathbb N$;
\item[$(5)$] $(X,\mathcal U(X))$ has the property $\bigcup_{fin}(\mathcal O,\Omega)$
($II\uparrow M(X,\mathcal U(X))$).
\end{itemize}
\end{thm}

In addition, the property $\bigcup_{fin}(\mathcal O,\Gamma)$ is a counterpart
of itself.

\begin{thm} \label{main-hur}
Let $X$ be a Tychonov space
and $\mathcal U$ be a natural uniformity on $A(X)$. Then the following conditions are equivalent:
\begin{itemize}
\item[$(1)$]  $(F(X),\mathcal U^\ast)$ has the property $\bigcup_{fin}(\mathcal O,\Gamma)$;
\item[$(2)$] $(F(X)^n,\mathcal U^\ast)$ has the property $\bigcup_{fin}(\mathcal O,\Gamma)$ for all $n\in\mathbb N$;
\item[$(3)$]  $(A(X),\mathcal U^\ast)$ has the property $\bigcup_{fin}(\mathcal O,\Gamma)$;
\item[$(4)$] $(A(X)^n,\mathcal U^\ast)$ has the property $\bigcup_{fin}(\mathcal O,\Gamma)$ for all $n\in\mathbb N$;
\item[$(5)$] $(X,\mathcal U(X))$ has the property $\bigcup_{fin}(\mathcal O,\Gamma)$.
\end{itemize}
\end{thm}

The following satement i crucial.

\begin{lem} \label{Lind} Let $X$ be a Lindel\"of space. Then $X$
has the property $\bigcup_{fin}(\mathcal O,\mathcal X)$ if and only
if so is the uniform space $(X,\mathcal U(X))$, where $\mathcal X\in\{\mathcal O,\Omega,\Gamma\}$.

In addition, the second player has a winning strategy in the Menger game on
$X$ if and only if he has a winning strategy in this game on $(X,\mathcal U(X))$.
\end{lem}

Combinning Lemma~\ref{Lind} with Theorems~\ref{main} and \ref{main-hur},
we obtain the following important corollary
\begin{cor} \label{M-equiv}
\begin{itemize}
\item[$(1)$]
The properties  $\bigcup_{fin}(\mathcal O,\mathcal X)$ as well as the property
$II\uparrow M(X)$ are
$M$-invariant within the class of Lindel\"of topological spaces, where $\mathcal X\in\{\Omega,\Gamma\}$.
\item[$(2)$] If $X$ has one of the above  properties  and $X^n$ is Lindel\"of
for all $n\in\omega$, then so is every $A$-equivalent to $X$ space $Y$.
\item[$(3)$] If $X$ is hereditarily Lindel\"of and has one of the above properties,
then so is every $M$-equivalent to $X$ space $Y$.
\end{itemize}
\end{cor}
We refer the reader to \cite{Tk00} for corresponding definitions.

As it was shown in \cite{Zd}, it is consistent with ZFC that each  topological space
with the property $\bigcup_{fin}(\mathcal O,\mathcal O)$ is $\bigcup_{fin}(\mathcal O,\Omega)$.
Therefore,  Corollary~\ref{M-equiv} implies that it is consistent that
the  property $\bigcup_{fin}(\mathcal O,\mathcal O)$ is  $M$-invariant within the class of Lindel\"of spaces.
\begin{prob} \label{pr1}
Is it consistent that the property $\bigcup_{fin}(\mathcal O,\mathcal O)$ is not $M$-invariant within the class
of Lindel\"of spaces?
\end{prob}

\nby{Lubmoyr Zdomsky}

\section{An equivalent of \SPMBul{} 2's \probmonth{}}
A family  $\mathcal F $ of infinite subsets of a countable set $C$
is said to be a \emph{semifilter}, if it is closed under taking
supersets and finite modifications of its elements.

The powerset $\mathcal P(C)$ admits a natural structure of a
compact space, and thus we can speak about topological properties
of its subsets such as semifilters.
 In \cite{Zd} the following characteriation of the Hurewicz property  was proven.
\begin{thm}
A Lindel\"of topological space $X$ has the Hurewicz property if
and only if for every countable large cover $u$ of $X$ the
smallest semifilter $\mathcal F$ on $u$ containing the family
$\{\{U\in u:x\in U\}:x\in X\}$ is Hurewicz.
\end{thm}
This characterization permits us to prove the subsequent
reformulation of the problem whether the Hurewicz property implies
the property $S_{fin}(\Gamma,\Omega)$.
\begin{thm} The following conditions are equivalent:
\begin{itemize}
\item[$(1)$] The Hurewicz property implies
$S_{fin}(\Gamma,\Omega)$; \item[$(2)$] For every Hurewicz
semifilter $\mathcal F$ on $\omega\times\omega$ such that
$\{n\}\times\omega\subset^\ast F$ for every $F\in\mathcal F$ there
exists a sequence $(K_n)_{n\in\omega}$ of finite subsets of
$\omega$ such that each element of the smallest filter containing
$\mathcal F$ meets $\bigcup_{n\in\omega}\{n\}\times K_n$.
\end{itemize}
\end{thm}

\nby{Lubmoyr Zdomsky}

\section{Boise Extravaganza In Set Theory (March 25--27, 2005)}
We are pleased to announce our fourteenth annual BEST conference.
There will be 4 talks by invited speakers:
Mirna Dzamonja (University of East Anglia, England);
Greg Hjorth (UCLA);
Paul Larson (Miami University);
William Mitchell (University of Florida).

The talks will be held on Friday, Saturday and Sunday in the
Department of Mathematics at Boise State University.
BEST social events are planned as well.

The conference webpage at
\url{http://math.boisestate.edu/~best/}
contains the most current information including lodging, abstract
submission, maps, schedule, etc.
Anyone interested in giving a talk and/or participating
should contact one of the organizers
(Justin Moore and Bernhard Koenig)
as soon as possible.
You can reach both of us by writing to
\texttt{best@math.boisestate.edu}

The conference is  supported by a grant from the National
Science Foundation, whose assistance is gratefully acknowledged.

\nby{Justin Moore}

\section{Problem of the Issue}

Recall that $\le^*$ is defined on $\NN$ by:
$f\le^* g$ if $f(n)\le g(n)$ for all but finitely many $n$.
Let $\b$ (the \emph{unbounding number}) denote the
minimal cardinality of an unbounded subset of $\NN$
(with respect to $\le^*$).
\begin{prob}
Is it provable that there is a set of reals of cardinality $\b$
which satisfies $\sone(\Gamma,\Gamma)$.
\end{prob}
The closest results we know of are the following.
The answer is positive if we replace $\b$ by $\t$
(Scheepers \cite{wqn}),\footnote{$\t$ is the \emph{tower number},
see Issue 5 of the \SPMBul{}.
$\t\le\b$, but it is consistent that $\t<\b$.}
or if we replace
$\sone(\Gamma,\Gamma)$ by ``$\ufin(\Gamma,\Gamma)$+no perfect subsets''
(Bartoszy\'nski-Tsaban \cite{ideals}).
It is still open whether ``$\ufin(\Gamma,\Gamma)$+no perfect subsets'' implies
(and is therefore equivalent to) $\sone(\Gamma,\Gamma)$ \cite{coc2}.

\section{Problems from earlier issues}
In this section we list the still open problems among
the past problems posed in the \SPMBul{}
(in the section \probmonth{}).
For definitions, motivation and related results, consult the
corresponding issue.

For conciseness, we make the convention that
all spaces in question are
zero-dimentional, separable metrizble spaces.

\begin{issue}
Is $\binom{\Omega}{\Gamma}=\binom{\Omega}{\Tau}$?
\end{issue}

\begin{issue}
Is $\ufin(\Gamma,\Omega)=\sfin(\Gamma,\Omega)$?
And if not, does $\ufin(\Gamma,\Gamma)$ imply
$\sfin(\Gamma,\Omega)$?
\end{issue}

\stepcounter{issue}

\begin{issue}
Does $\sone(\Omega,\Tau)$ imply $\ufin(\Gamma,\Gamma)$?
\end{issue}

\begin{issue}
Is $\p=\p^*$? (See the definition of $\p^*$ in that issue.)
\end{issue}

\begin{issue}
Does there exist (in ZFC) an uncountable set satisfying $\sone(\BG,\B)$?
\end{issue}

\stepcounter{issue}

\begin{issue}
Does $X \nin \NON(\M)$ and $Y\nin\mathsf{D}$ imply that
$X\cup Y\nin \COF(\M)$?
\end{issue}

\begin{issue}
Is $\split(\Lambda,\Lambda)$ preserved under taking finite unions?
\end{issue}
\begin{proof}[Partial solution]
Consistently yes (Zdomsky). Is it ``No'' under CH?
\end{proof}

\begin{issue}
Is $\cov(\M)=\fo$? (See the definition of $\fo$ in that issue.)
\end{issue}

\begin{thebibliography}{00}
\bibitem{ideals}
T.\ Bartoszy\'nski and B.\ Tsaban,
\emph{Hereditary topological diagonalizations and the Menger-Hurewicz Conjectures},
Proceedings of the American Mathematical Society, to appear.
\arx{math.LO/0208224}

\bibitem{En}
R.\ Engelking,
\emph{General Topology. Revised and completed edition},
Sigma Series in Pure Mathematics \textbf{6}(1989), Heldermann
Verlag, Berlin.

\bibitem{HRT}
C.\ Hernandes, D.\ Robbie, M.\ Tkachenko, \emph{some properties of $o$-bounded and strictly $o$-bounded
topological groups}, Appl.\ General Topology \textbf{1}(2000),
29-43.

\bibitem{coc2}
W.\ Just, A.\ W.\ Miller, M.\ Scheepers, and P.\ J.\ Szeptycki,
\emph{The combinatorics of open covers II},
Topology and its Applications \textbf{73} (1996),
241--266.

\bibitem{Sch}
M.\ Scheepers, \emph{Combinatorics of open covers I: Ramsey Theory},
Topology Appl. \textbf{69} (1996) 31-62.

\bibitem{wqn}
M.\ Scheepers,
\emph{Sequential convergence in ${\sf C}_p(X)$ and a covering property},
East-West Journal of Mathematics \textbf{1} (1999),
207--214.

\bibitem{Tk98}
M.\ Tkachenko, \emph{Introduction to topological groups}, Topology
and Appl. \textbf{86}(1998), 179-231.

\bibitem{Tk00}
M.\ Tkachenko, \emph{Topological groups for topologists: part II},
Bol. Soc. Mat. Mexicana (3) Vol. 6, 2000.

\bibitem{Zd}
L.\ Zdomsky,
\emph{A semifilter approach to selection principles},
submitted.

\ed